\documentclass[12pt]{amsart}
\usepackage{amssymb}
\usepackage{latexsym}
\usepackage{epsf}
\usepackage{amssymb, epic,eepic,epsfig,amsbsy,amsmath,amscd}
\newtheorem{thm}{Theorem}[section]
\newtheorem{dfn}[thm]{Definition}
\newtheorem{Thm}{Theorem}
\newtheorem{pro}[thm]{Proposition}
\newtheorem{lem}[thm]{Lemma}

\theoremstyle{definition}
\newtheorem{definition}{Definition}[section]

\def\1{{\rm1\mathchoice{\kern-0.25em}{\kern-0.25em}
        {\kern-0.2em}{\kern-0.2em}I}}

\newcommand{\lmn}[1]{\vadjust{\setbox1=\vtop{\hsize 25mm
\parindent=0pt\baselineskip=9pt
\rightskip=4mm plus 4mm#1}
\hbox{\kern-26mm\smash{\raise .5ex\box1}}}}
\input epsf.sty

\newcommand{\nc}{\newcommand}
\def\be#1\ee{\begin{equation}#1\end{equation}}
\nc{\bc}{\begin{center}} \nc{\ec}{\end{center}} \nc{\bb}{\mathbb}
\nc{\cal}{\mathcal} \nc{\frk}{\mathfrak} \nc{\N}{{\mathsf N}}
\nc{\K}{{\mathsf K}} \nc{\fk}{\mathbf{k}} \nc{\fn}{\mathbf{n}}
\nc{\fb}{\mathbf{b}}  \nc{\e}{\varepsilon} \nc{\ev}{{\rm{ev}}}

\theoremstyle{remark}

\def\Z{{\mathbb Z}}

\def\N{{\mathbb N}}

\def\v8{\vskip 8pt}

\def\x{\mathbf x}
\def\y{\mathbf y}
\def\z{\mathbf z}

\def\S{\mathbf S}
\def\Rect{{\text{Rect}^0}}

\def\Field{\mathbb F}
\newcommand\Ta{\mathbb{T}_\alpha}
\newcommand\Tb{\mathbb{T}_\beta}

\renewcommand\th{^{\text{th}}}
\def\HFK{\widehat{\operatorname{HFK}}}
\def\HFL{\widehat{\operatorname{HFL}}}

\def\lo{_{\text{long}}}
\def\sh{_{\text{short}}}

\def\BR{\mathbb R}

\newcommand\Xs{\mathbb{X}}
\newcommand\Os{\mathbb{O}}

\def\int{\mathrm {int}}
\def\NESW{\mathcal I}

\def\Field{\mathbb F}

\newcommand\alphas{\mbox{\boldmath$\alpha$}}
\newcommand\betas{\mbox{\boldmath$\beta$}}
\newcommand\ws{\mathbf w}
\newcommand\zs{\mathbf z}

\def\H{\mathbb H}

\def\cm{\cdot}
\def\x{\mathbf x}
\def\y{\mathbf y}

\def\Bi{{\text{Bigon}^0}}

\begin{document}

\title[Combinatorial Link Floer Homology]
{A Simplification of  Combinatorial \\[.1cm]
Link Floer Homology}
\author{Anna Beliakova}
\address{Institut f\"ur Mathematik, Universit\"at Z\"urich,
 Winterthurerstrasse 190,
CH-8057 Z\"urich, Switzerland}
\email{anna@math.unizh.ch}

\email{}


\begin{abstract}
We define a new combinatorial complex computing the hat version of
 link Floer homology over $\Z/2\Z$, which turns out to be
   significantly  smaller
than the  Manolescu--Ozsv\'ath--Sarkar one.
\end{abstract}

\maketitle
\section*{Introduction}
Knot Floer homology is a powerful knot invariant  constructed
by Ozsv\'ath--Szab\'o \cite{Knots} and Rasmussen \cite{RasmussenThesis}.
In its basic form,  the knot Floer homology 
$\HFK(K)$ of a knot $K\subset S^3$ 
is a finite--dimensional bigraded vector space
over $\Field=\Z/2\Z$ 
$$\HFK(K)=\;\bigoplus\limits_{d\in \Z,i\in \Z} \;\HFK_d (K,i)\, ,$$
where $d$ is the  Maslov and $i$ is the Alexander grading.
Its graded Euler characteristic 
$$\sum_{d,i} (-1)^d \mathrm {rank}\;\HFK_d(K,i)t^i=\Delta_K(t)$$
is equal to the symmetrized Alexander polynomial  $\Delta_K(t)$.
The knot Floer homology enjoys the following
symmetry extending that of the Alexander polynomial.
\be\label{sym}
\HFK_d(K,i)=\HFK_{d-2i}(K,-i)\ee

By the result of Ozsv\'ath--Szab\'o \cite{GenusBounds},
the maximal Alexander grading $i$, such that $\HFK_*(K,i)\neq 0$ is 
 the Seifert genus $g(K)$ of $K$. Moreover, 
Ghiggini showed for $g(K)=1$ \cite{Ghiggini} and
Yi Ni in general \cite{YiNi},
that
the knot is fibered if and only if $\mathrm {rank}\, \HFK_*(K,g(K))=1$. 
A concordance invariant bounding from below the slice genus of the knot
can also be extracted from  knot Floer homology \cite{4BallGenus}.
For torus knots the bound is sharp, providing a new proof of the  Milnor
conjecture. The first  proof of the Milnor conjecture
was given by Kronheimer and Mrowka  \cite{KMMilnor},
then Rasmussen  \cite{RasmussenSlice}  proved it combinatorially 
by using
  Khovanov homology \cite{Kho}.

Knot Floer homology was extended to links in \cite{Links}.
The first combinatorial construction of the link Floer homology
was given in \cite{MOS} over $\Field$ and then  in  \cite{MOST}
 over $\Z$. Both constructions use grid diagrams of links.

A {\em  grid diagram}  is a square
grid on  the plane with $n \times n$ squares. Each square is decorated
either with an $X$, an $O$, or nothing. Moreover, 
every row and every column contains exactly one $X$ and one $O$.
The number $n$ is called  {\em complexity} of the diagram.
Following \cite{MOST}, we denote
the set of all $O$'s and $X$'s by~$\Os$ and~$\Xs$, respectively.

\begin{figure}
\mbox{\epsfysize=6cm \epsffile{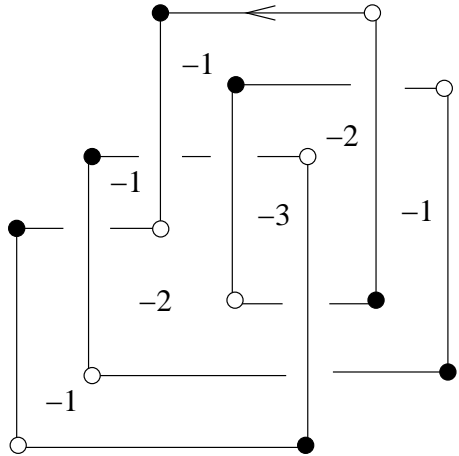}} 
\v8

\caption{{\bf A rectangular diagram  for $5_2$ knot.} The number
associated to a domain is  minus the winding 
number for its points. The sets $\Xs$ and $\Os$ consist of black and
white points, respectively.}
\label{diag}
\end{figure}

Given a grid diagram, we
construct an oriented, planar link projection  by drawing horizontal
segments from the $O$'s to the $X$'s in each row, and vertical
segments from the $X$'s to the $O$'s in each column. We assume that at every
intersection point  the vertical  segment  overpasses 
the horizontal one. This produces a planar {\em
rectangular diagram}  $D$ for an
oriented link $L$ in $S^3$. Any link in $S^3$ admits a rectangular 
diagram (see e.g. \cite{Dy}). An example is shown in Figure
\ref{diag}.

In \cite{MOS}, \cite{MOST} the grid lies on the torus, 
obtained by gluing the top most segment of the grid to the bottommost one
and the leftmost segment to the right most one.
In the torus, the horizontal
and vertical segments of the grid become circles.
The  MOS complex is then generated by $n$--tuples of 
intersection points between
 horizontal and vertical circles, such that exactly one point
belongs to each horizontal (or vertical) circle. 
The differential is defined as follows:
$$
\partial \x =\sum_{\y\in S_\x} \;\;\;\sum_{r  \in \Rect(\x,\y)} \;\;\y\, ,
$$
where $S_\x$ is the subset of generators 
that have $n-2$ points in common with $\x$. For $\y\in S_{\x}$, 
$\Rect(\x,\y)$ is the set of rectangles 
 with vertices 
 $\x\setminus(\x\cap\y)$ and  $\y\setminus(\y\cap\x)$, 
  whose interior
 does not contain 
 $X$'s and $O$'s or points among $\x$ and $\y$. Moreover,  
a counterclockwise
rotation along the arc of the horizontal oval, leads from the vertices 
in $\x$ to the ones in $\y$.

The Alexander grading is given by formula \eqref{eq:Alexander} below, 
and the Maslov grading by \eqref{mas} plus one.
The MOS complex has $n!$ generators. This number  greatly exceeds the rank 
of its  homology. For the trefoil, for example,
the number of generators is 120, while the rank of
$\HFK(3_1)$ is 3.

In this paper, we construct another combinatorial
complex computing link Floer homology, which has significantly
less generators.
All knots with less than 6 crossings admit rectangular diagrams
where  all differentials
 in our complex are
zero, and the rank of the homology group is equal to the number of 
 generators.

\subsection*{Main results}
\begin{figure}
\mbox{\epsfysize=6.5cm \epsffile{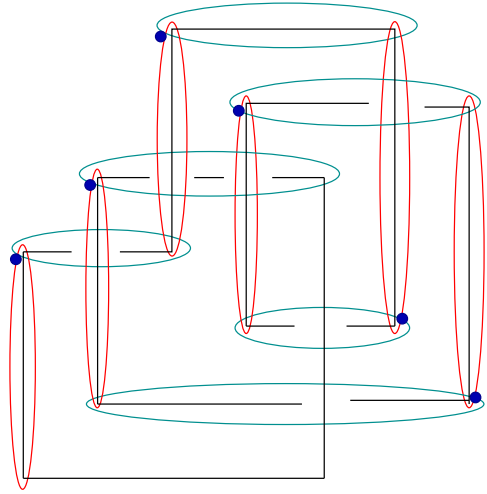}}
\v8

\caption{{\bf Collection of short ovals  for $5_2$ knot.} 
The dots show a generator
in Alexander grading 1.}
\label{fig:5-2}
\end{figure}

Our construction also uses rectangular diagrams.
Given an oriented link $L$ in $S^3$,
let $D$ be its  rectangular diagram in $\BR^2$. Let
us  draw $2n-2$
 narrow short ovals around all but one horizontal and all but one
vertical segments of the rectangular  diagram $D$ in such a way,
that the outside domain has at least one point among $\Xs$ or $\Os$.
We denote by $\S$ the set of unordered
 $(n-1)$--tuples of intersection points between the horizontal
and vertical  ovals, such that exactly one point belongs to each
horizontal (or vertical) oval. We assume throughout this paper
that the ovals intersect transversely.
An example is shown in Figure \ref{fig:5-2}.

A chain complex $(C\sh(D),\partial\sh)$
computing the hat version of 
 link Floer homology of $L$ over  $\Field=\Z/2\Z$
is  defined as follows.
The generators are elements of $\S$.
The bigrading on $\S$ can be constructed analogously to those
in \cite{MOS}.
Suppose   $\ell$ is the number of components of $L$.
Then the Alexander grading 
is a function $A \colon \S \longrightarrow (\frac{1}{2}\Z)^\ell$, 
defined as follows.

First, we define a function $a:\S \to 
\Z^{\ell}$. For a  point $p$, the $i\th$ component
of $a$ is minus  the winding number of the projection of the
$i\th$ component of the 
oriented link around $p$. In the grid diagram, we 
have $2n$ distinguished squares  containing $X$'s or $O$'s. Let
$\{c_{i,j}\}, \ i\in\{1,\ldots,2n\},\ j\in\{1,\ldots,4\}$, be the
vertices of these squares.
 Given $\x\in\S$, we set
\begin {equation}
\label{eq:Alexander}
A(\x)=\sum_{x\in \x} a(x)- \frac{1}{8} \Bigl(\sum_{i,j} 
a(c_{i,j})\Bigr) - 
\left(\frac{n_1-1}{2},\ldots,\frac{n_\ell-1}{2}\right),
\end {equation}
where here $n_i$ is the complexity of the $i\th$~component of~$L$, i.e.
 the number of horizontal segments belonging to this component.

The homological or Maslov  grading is  a function $M:\S\to\Z$
 defined as follows. Given two collections $A$,
$B$ of finitely many points in the plane, let $\NESW(A,B)$ be the
number of pairs $(a_1,a_2)\in A$ and $(b_1,b_2)\in B$ with $a_1<b_1$
and $a_2<b_2$. Let $J(A,B):=1/2(\NESW(A,B)+ \NESW(B,A))$.  
Define 
\be\label{mas}
M(\x)=J(\x,\x)-2J(\x,\Os)
+J(\Os,\Os).\ee


\begin{figure}
\mbox{\epsfysize=8.5cm \epsffile{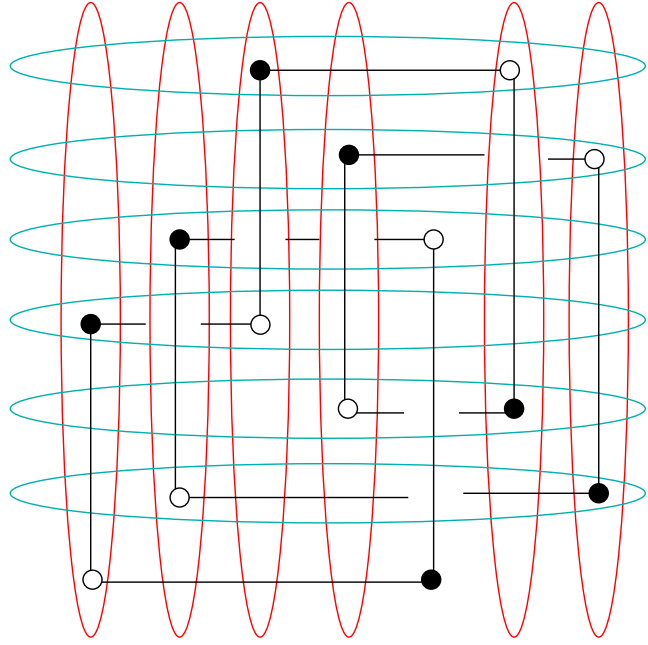}} 
\v8
\caption{\bf Collection of long ovals.}
\label{long}
\end{figure}

 
To construct  a differential $\partial\sh$
we first need
to consider the 
 complex $(C\lo(D),\partial\lo)$ defined in the same way
as $(C\sh(D),\partial\sh)$
but where the ovals are  as long as $n\times n$ grid.
An example is shown in  Figure \ref{long}.
The differential
$$\partial\lo(\x)=\sum_{\y\in S_\x\cup S'_\x}\;\;\;\;
\sum_{r\in \Rect(\x,\y)\cup\Bi(\x,\y)}
\;\;\;\y 
$$
where $S'_\x$ is the subset of  generators that have
 $n-2$ points in common with $\x$ and for $\y\in S'_\x$,  
$\Bi(\x,\y)$ is the  bigon with
 vertices $\x\setminus(\x\cap\y)$ and $\y\setminus(\y\cap\x)$,
whose interior does not contain $X$'s  and $O$'s. 
Moreover,  a counterclockwise
rotation along the arc of the horizontal oval, leads from the vertex 
in $\x$ to the one in $\y$.

Note that $(C\lo(D),\partial\lo)$ coincides with
the   complex  $(C(S^2,\alphas,\betas, \Xs,\Os),\partial)$
 defined by Ozsv\'ath--Szab\'o
in \cite{Links}, where the vertical and horizontal ovals are identified 
with $\alphas$
and $\betas$ curves, respectively, and  $\Xs$ and $\Os$ are  extra basepoints.
Like the MOS complex, this complex is combinatorial, since
all domains suitable for the differential are either rectangles or bigons
(compare \cite{SW}).
Following \cite{Links}, 
we will call elements of $\Xs\cup\Os$ {\it basepoints}
in what follows.


To this complex 
we further apply  a simple lemma from homological
algebra, that allows us to construct a homotopy equivalent complex
$(C\sh(D),\partial\sh)$
with $\S$ as a set of generators
(or to shorten the  ovals $\alphas$ and $\betas$ keeping 
track of the differential).
 
In Section 1, we describe an algorithm that for
any domain connecting two generators decides
 whether 
it counts for the differential $\partial\sh$ or not. 
Furthermore, we distinguish
 a large class
of domains that always count.
 In general, however
the count   depends on the order of shortening of ovals, which
replace  the choice of a complex structure in the
analytic setting.



Let $V_i$ be the two dimensional  bigraded vector space 
over $\Field$ spanned  by one
generator in Alexander and Maslov gradings zero and
another one in Maslov grading $-1$ and Alexander grading minus the $i$--th
basis vector.

\parbox[t]{16cm}{
\begin{Thm}\label{main}
Suppose $D$ is a rectangular diagram of an oriented $\ell$ component link $L$, 
where  the $i^{th}$ component of
  $L$ has complexity  $n_i$.  
Then the homology
$$H_*(C\sh(D),\partial\sh)=\HFL(L)\otimes 
\bigotimes\limits^\ell_{i=1}V_i^{n_i-1}$$ 
can be computed algorithmically.
\end{Thm}}


The complex $(C\sh(D),\partial\sh)$
 has much fewer
 generators than  the MOS complex (compare Section 3).
Recently, Droz
introduced signs in our construction \cite{D1} and wrote  
a computer program   realizing
our algorithm over $\Z$  \cite{Droz}. His program
allows to determine the Seifert genus and fiberedness
of  knots   until  16 crossings and also 
 to study the torsion part of knot Floer homology.

The paper is organized as follows.
 Theorem \ref{main} is proved in Section \ref{proofs}.
In Section \ref{ind}
 we  introduce a big class of domains that always count for the
differential.
 In the last section  we compute Floer homology of $5_2$ knot and 
discuss further computations made by  Droz.

\subsection*{Acknowledgments} First of all, I wish to thank Stephen 
Bigelow for explaining to  me his homological constructions of
link invariants, which provided  the original motivation for this project,
and also for  regular conversations and helpful hints.

 I would  like
to express my gratitude to Jean--Marie
Droz  for many interesting discussions and, especially, for 
writing a program computing the homology of the complex
constructed in this paper.  

I profited a lot from discussions with Stephan Wehrli.
Special thanks go to
Ciprian Manolescu for sharing his knowledge of Heegaard Floer homology
and for his valuable suggestions after reading the preliminary version
of this paper.


\section{The complex 
 $(C\sh(D),\partial\sh)$}\label{proofs}

\subsection{Intermediate complex $(C,\partial)$} Suppose $D$ is a rectangular
diagram  of complexity $n$  for an oriented link $L$.
Let 
$(C,\partial)$ be the complex
generated over $\Field$ by $(n-1)$--tuples of intersection points
 between horizontal and vertical ovals, such that exactly one
 point belongs to each horizontal (or vertical) oval
as defined  in the Introduction. 
The length of the ovals can be intermediate between long and short ones.
We also assume that the outside domain has
at least one basepoint inside.


Given $\x,\y\in S(C)$, there is an oriented closed curve
$\gamma_{\x,\y}$ composed of arcs belonging to 
horizontal and vertical ovals, where each piece of a  horizontal oval 
connects a point in $\x$ to a point in $\y$ (and hence each piece of the 
vertical one goes from a point in $\y$ to a point in $\x$).
In $S^2$, there exists an oriented (immersed)
 domain $D_{\x,\y}$
 bounded by $\gamma_{\x,\y}$. 
The points in $\x$ and $\y$ are called
{\em corners} of  $D_{\x,\y}$.

 Let $D_i$ be  the closures 
of the connected components
of the complement of ovals in $S^2$. Suppose that the orientation of $D_i$ 
is induced by the orientation of $S^2$.
Then we say that a domain  $D=\sum_i {n_i} D_i$ connects two generators
if for all $i$,  $n_i\geq 0$ and $D$ is connected. 
Let $\mathfrak D$ be the set of all  domains connecting two generators
which  contain neither corners nor points among $\Xs$ and $\Os$  inside.

We define 
\be\label{different}
\partial \x :=\sum_{M(\y)=M(\x)-1} \;\sum_{D_{\x,\y}\in \mathfrak D}\; 
m(D_{\x,\y})\; \y\, ,\ee
where $M(\x)$ is the Maslov grading defined by (\ref{mas})
and
$m(D_{\x,\y})\in \{0,1\}$ is a  multiplicity of $D_{\x,\y}$. 
We set
$m(D_{\x,\y})=1$ for rectangles and bigons without basepoints inside.
This defines the differential for the long oval complex.
In general, $m(D_{\x,\y})$ can be defined  by using
the procedure of  shortening of ovals described in the next section
inductively. In particular,  $\partial\sh$
 is defined in Section 1.3.
\vspace{.2cm}

\subsection{Shortening of ovals}
\begin{figure}
\mbox{\epsfysize=1.8cm \epsffile{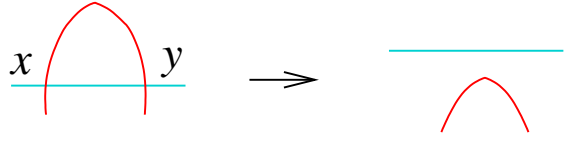}} 
\v8

\caption{\bf Removing of a bigon without basepoints inside.}
\label{fig:bigon}
\end{figure}

Assume that one vertical and one horizontal oval used to define
$(C,\partial)$ form a bigon 
with corners $x$ and $y$  without basepoints inside. We further assume that
a counterclockwise rotation along the arc of the horizontal oval, leads 
from $x$ to $y$.
An example is shown  on the left of Figure \ref{fig:bigon}.
Suppose that $\partial$ is given by \eqref{different} with known multiplicities.

Let $(C',\partial')$ be a new complex obtained from $(C,\partial)$
as follows. The set of generators $S(C')$ is obtained from $S(C)$
by removing all generators containing $x$ or $y$.
The differential 
 $$\partial':=P\circ(\partial +\partial \circ h\circ\partial)\circ I\, ,$$
where $I: C'\to C$ and $P:C\to C'$  are
the obvious inclusion and projection. Moreover, for any $\x \in S(C)$, 
$h(\x)$ is zero whenever $y\not\in \x$, otherwise $h(\x)$ is obtained from $\x$
by  replacing $y$ by $x$.

\begin{pro}\label{bigon}
 $(C',\partial')$ is a chain complex homotopic to
$(C,\partial)$. The differential $\partial'$ can be 
given by \eqref{different} for a new system of curves obtained by
shortening of one oval as shown in  Figure \ref{fig:bigon}.
\end{pro}

For the  proof of this proposition we will need the following lemma.

\begin{lem}\label{sublemma}
Whenever $x\in \x$,
we have $h\circ \partial \x=\x$ 
\end{lem}

\begin{proof}
We have to show that
for $\x$ with $x\in \x$,  $\z$ with $y\in \z$
occurs in $\partial \x$
 if and only if
  $\z$ is obtained from $\x$ by replacing
$x$ by $y$.  

Indeed, if $\z\in \partial\x$, then there exists a domain
$D_{\x,\z}\in \mathfrak D$  with $m(D_{\x,\z})=1$.
Let us construct its boundary $\gamma_{\x,\z}=\partial D_{\x,\z}$.
We start at $x$ and go to $y$ along an arc of a horizontal oval.
No other points in $\x$ or $\z$ 
belong to this arc, since only one point of each generator
is
 on  the same oval. There are two choices for this arc: 
the short or the long one. The  domain having the long arc as 
a part of its boundary
 contains basepoints inside and will not be counted.
Analogously, going back from $y$ to $x$ we have to take the short arc
of the vertical oval, otherwise the domain will not be in $\mathfrak D$.
Hence,  $D_{\x,\z}$ contains a small bigon between $x$ and $y$ 
as a boundary component.
Let us show that  $D_{\x,\z}$ has only one boundary component.
Clearly, there are no further boundary components inside of this bigon.
On the other hand, the bigon could not be an inner boundary component,
hence in this case, its orientation would be reversed, and then
 $D_{\z,\x}$ would be in $\mathfrak D$, but not
 $D_{\x,\z}$. Since $\mathfrak D$  contains connected domains only,  
we proved the claim.

\end{proof}

\subsection*{Proof of Proposition \ref{bigon}}

Let us define
the  maps $F: C\to C'$  and $G:C' \to C$
 as follows.
 $$F=P\circ (\1 +\partial\circ h)\;\;\;\; G=(\1 +h\circ\partial)\circ I\, $$ 
Here $\1$ is the identity map. It is not difficult to see, that
 $F \circ G$ is the identity map on $C'$ since
$P\circ h$, $h\circ I$ and $h^2$ are zero maps.




The map
$G\circ F$ is homotopic to the identity on $C$, i.e.
$
G\circ F + \1=\partial\circ  h+h\circ \partial$.
 This is easy to see  for generators in
$S(C)\cap S(C')$ or for $\x$ with $x\in \x$. 
Let us assume $y\in \y$, then $\partial\circ h(\y)=\y+\tilde\x+\tilde \z$,
where $\tilde \x$ is a linear combination of generators, such that 
any of them contains $x$ and $\tilde\z$ is a linear combination
of generators in $S(C)\cap S(C')$. Then 
$$G\circ F (\y)=(\1+h\circ\partial)(\tilde \z)= \tilde \z +h(\partial \y)
+\tilde \x$$
by using $\partial^2=0$ and the previous claim. On the other hand,
 $$(\1 +h\circ \partial +\partial\circ h)(\y)=
\y+\y+\tilde\x+\tilde\z+h(\partial \y)\, .$$

Now using Lemma \ref{sublemma} one can easily show that
the differential
 $\partial'$ coincides with $F \circ\partial \circ G$.
Furthermore, using the homotopy to the identity, proved above, we derive
 $\partial'^2= F \circ\partial \circ G \circ F \circ\partial \circ G =0$. 
This shows that $(C', \partial')$ is indeed a 
chain complex.

To show that $F$ and $G$ are chain maps, i.e.
 $\partial'\circ F= F\circ\partial$ and $
\partial\circ  G=G\circ\partial'$,  
we again use  $\partial \circ h\circ\partial\circ h\circ \partial =
\partial \circ h\circ\partial$, which is a consequence of
Lemma \ref{sublemma}.


Finally, we would like  to show
 that the new differential
 $\partial'=P\circ(\partial +\partial\circ h\circ\partial)\circ I$ 
can be realized by counting of Maslov index one domains for a new 
system of curves.
We write $\y\in \partial\x$ if $m(D_{\x,\y})=1$. 
Assume $\y\in \partial\x$ and $\x$, $\y$ do not
contain the corners of the bigon $x$ and $y$. 
 For all such $\x$ and $\y$, 
we also have $\y\in\partial'\x$, and they are connected by a Maslov index 
one domain $D_{\x,\y} \in \frak D$.

Furthermore,
assume $\mathbf a,\mathbf b \in \S(C)$ do
not contain $x$ and $y$, then for any $\x,\y\in \S(C)$ with $y\in \y$, 
$x\in \x$, such that
$h(\y)=\x$,
 $\y\in\partial\mathbf a$,
 and $\mathbf b\in\partial\x$, 
$\mathbf b$ occurs once in $\partial'\mathbf a$. 
Note that $\mathbf b \not\in \partial \mathbf a$, since
any domain  connecting $\mathbf a$ to $\mathbf b$
contains either the bigon with negative orientation or basepoints,
hence they do not count for the differential.
However, the new system of ovals contains a domain 
 connecting $\mathbf a$ to 
$\mathbf b$ which is obtained from $D_{\mathbf a,\y}\cup D_{\x,\mathbf b}$ by
shortening the oval. 

\begin{figure}
\mbox{\epsfysize=5cm \epsffile{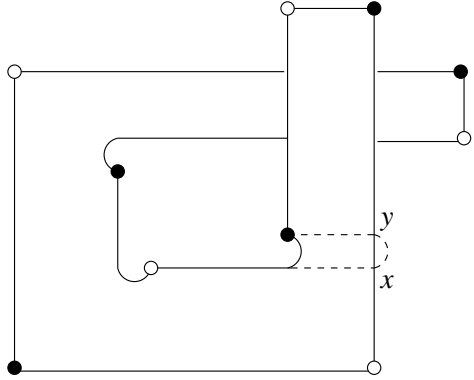}} 
\v8

\caption{{\bf Immersed polygon realizing a differential
from black to white points.} The prolongated oval is 
shown by a dashed line.}
\label{poly}
\end{figure}

This process is illustrated in Figure \ref{poly}, where
$\mathbf a$ and $\mathbf b$ are given by black and white points respectively;
$\y$ is obtained from $\mathbf a$ 
by switching the black point on the dashed oval 
to $y$ and the upper black point to the white point  on the same oval;
$\x$ is obtained from $\y$ by switching $y$ to $x$.

\begin{figure}
\mbox{\epsfysize=4cm \epsffile{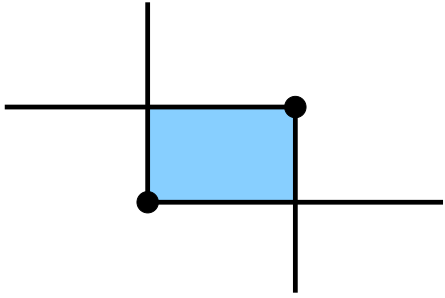}} 

\caption{{\bf Intersection of two domains with two corners on its boundary.} }
\label{intersect}
\end{figure}

It remains to show that the domain, obtained from $D_{\mathbf a,\y}\cup D_{\x,\mathbf b}$ by
shortening the oval,
belongs to  $\mathfrak D$.
The domain is obtained by connecting $D_{\mathbf a,\y}$ to $ D_{\x,\mathbf b}$ 
by two arcs shown on the right of Figure 4, i.e. the domain is connected.
It has no  basepoints inside, since
$D_{\mathbf a,\y}$ and $ D_{\x,\mathbf b}$  do not have them.
Analogously, it can be written as $\sum_i n_i D_i$ with all $n_i\geq 0$.
Finally, we show that our domain has no corners inside.
If  $D_{\mathbf a,\y}\cap  D_{\x,\mathbf b}$ is  empty, it follows 
from the assumption that $ D_{\x,\mathbf b},  D_{\mathbf a,\mathbf y}
\in \mathfrak D$, i.e. have no corners inside.
If the intersection is not empty, then its boundary  either contains
no corners or
at least two corners, one of them in $\y\setminus y=
\x\setminus x$  (see Figure \ref{intersect}).  
The last is impossible since 
$ D_{\x,\mathbf b}\in \mathfrak D$ has no corners of $\x$ 
and $D_{\mathbf a, \y}$ no corners of $\y$ inside.

\qed


{\bf Note.} The first statement of
Proposition \ref{bigon} is a particular case of the Gaussian 
elimination considered in e.g. \cite[Lemma 4.2]{Bar}.

\subsection{Definition $(C\sh(D),\partial\sh)$}
The complex $(C\sh(D), \partial\sh)$  is  obtained from
the complex with long ovals  by applying Lemma \ref{bigon} several times,
until the complex has $\S$ as the set of generators.
This subsection aims to give a recursive definition of the multiplicity
$m(D_{\mathbf a,\mathbf b})$ of $\mathbf b$ in the
differential $\partial\sh (\mathbf a)$ for all $\mathbf a\in \S$.  
  In general, it will depend 
 on the order   in which the ovals were shortened. 


Let us fix this order. If $D_{\mathbf a,\mathbf b}\in \mathfrak D$ 
is a polygon we count it with multiplicity one.
For any other domain in $\mathfrak D$  we prolongate the 
last oval that was
shortened to obtain this domain, and show whether in the resulting complex
$(C'(D),\partial')$,
one can find $\x$ and $\y$ with $\y\in \partial' \mathbf a$, $h(\y)=\x$ and
$\mathbf b\in \partial'\x$
as in the proof of Lemma
\ref{bigon}. If there is an odd number of such $\x$ and $\y$, then
$m(D_{\mathbf a,\mathbf b})=1$, otherwise the multiplicity is zero.
To determine  $m(D_{\mathbf a,\y})$ and $m(D_{\x,\mathbf b})$ 
in the new complex, we
 prolongate the next oval, and continue to do so until
the domains in question are polygons which always count.




In fact, the algorithm terminates already when the domains 
in question are
strongly indecomposable  without bad components  (as defined
in the next section) since they all count for the differential (cf. Theorem 
2.1 below).



\subsection{Proof of Theorem \ref{main}}
It remains to compute the homology of our complex.

By Lemma \ref{bigon},
 the complex $(C\sh(D),\partial\sh)$ is homotopy equivalent
to the complex with long ovals. The homotopy preserves
both gradings, since $\partial\sh$ count only Maslov index one domains
without basepoints inside.
The complex $(C\lo(D),\partial\lo)$ coincides with
the complex computing the hat version of  link Floer homology
from the genus zero Heegaard splitting of $S^3$
with extra basepoints (see \cite{MOS} for summary).
The relative Maslov and Alexander of these two complexes are also the same.
Moreover, the absolute
Alexander and Maslov gradings in our complex are fixed in such a way, that
$$\chi(C\sh(D),\partial\sh)=\left\{ \begin{array}{ll}
 \prod^{\ell}_{i=1} 
 (1-t^{-1})^{n_i-1}\Delta_L(t_1,...,t_\ell)& \ell>1 \\
 (1-t^{-1})^{n_i-1}\Delta_{L}(t) & \ell=1 \end{array}
 \right. $$
This can be shown by comparing our and MOS complexes 
(see \cite[Theorem 4.2]{D1}
for more details). 

Hence  Proposition 2.4 in \cite{MOS} computes the homology of our 
combinatorial complex.


\qed

\section{Domains that  always count}\label{ind}
To run the algorithm defining $\partial\sh$ we have to decide
at each step which domains count and which do not. 
In this Section we simplify the algorithm
by selecting a large class of domains 
that always count for the differentials obtained
from $\partial\lo$ by applying Lemma \ref{bigon} several times.
Let us start with some definitions.

\subsection{Maslov index}
 Let $e(S)$ be the Euler measure of 
a surface  $S$, which  for any surface $S$ with
$k$   acute right--angled  corners,
$l$  obtuse  ones, and Euler characteristic $\chi(S)$ is equal to 
$\chi(S) -k/4 +l/4$. Moreover, the Euler measure is additive under disjoint 
union and gluing along boundaries.
In \cite[Section 4]{Lip}, Lipshitz gave a formula computing
the  Maslov index 
 $M(D_{\x,\y})$ of $D_{\x,\y}$ as follows. 
\be\label{mas-domain}
M(D_{\x,\y})=e(D_{\x,\y}) +n_{\x}+n_{\y}\, ,
\ee
where
 $n_{\x}=\sum_{x\in \x} n_{x}$.
The number $n_{x}$ is the local multiplicity of the domain at the corner 
$x$, e.g.
 $n_{x}=0$
for an isolated corner, $n_{x}=1/4$
for an acute (or $\pi/2$--angled) corner, 
$n_{x}=1/2$ for a straight (or $\pi$--angled) corner  or $n_{x}=3/4$ 
for an obtuse (or
 $3\pi/2$--angled)
 one.
For  a composition of two domains
 $D_{\x,\z}=D_{\x,\y}\circ D_{\y,\z}$,
we have $M(D_{\x,\z})=M(D_{\x,\y})+ M(D_{\y,\z})$.

A path in a domain starting at an obtuse  or straight corner and
following a horizontal or vertical oval until the boundary of the domain
will be called a {\it cut}. There are two cuts
at any obtuse corner  and  one at any straight corner.

\begin{figure}
\mbox{\epsfysize=5cm \epsffile{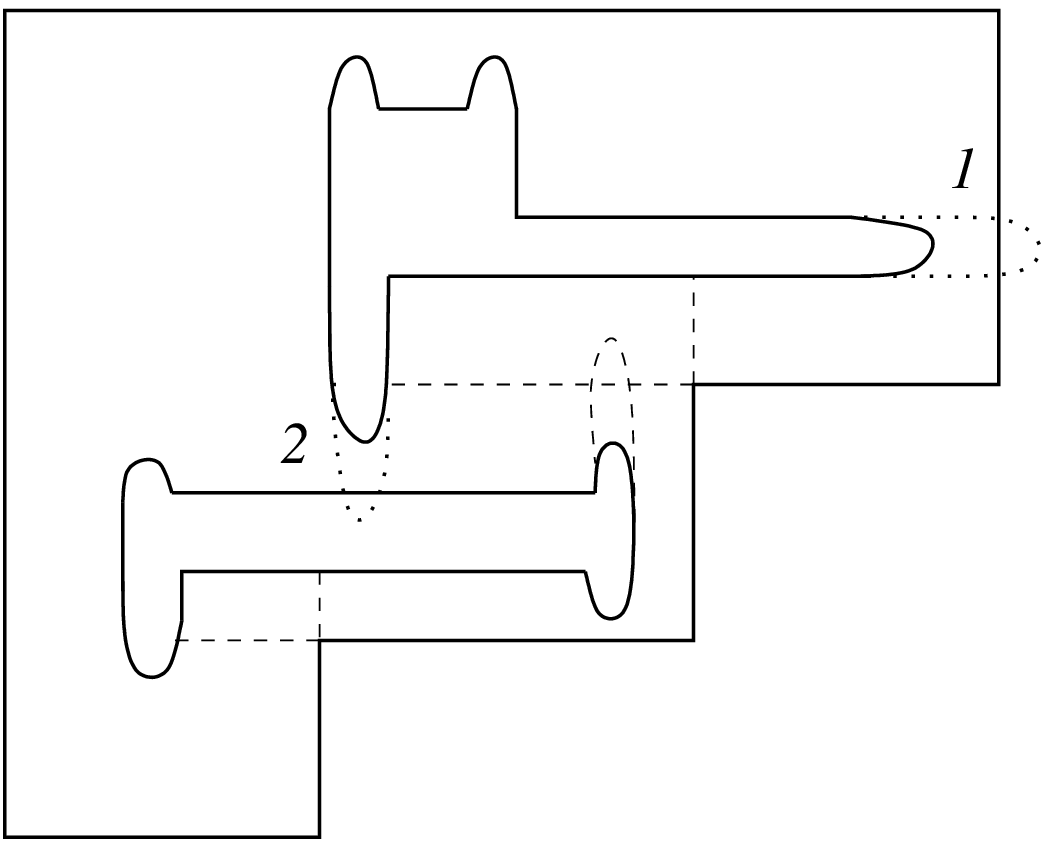}} 
\v8

\caption{\bf  Indecomposable, but not strongly indecomposable domain.}
The oval destroying the indecomposability is shown by the dashed line
without any number.
\label{notstrong}
\end{figure}

A  domain $D$ is called {\it decomposable} if it is a composition 
of Maslov index
zero and one domains; any other Maslov index 1 domain is called 
indecomposable.
A domain is called {\it strongly indecomposable} if 
the following conditions are satisfied:
\begin{itemize}
\item
it is indecomposable;
\item
no prolongations 
of ovals inside this domain destroy its indecomposability;
 \item
the cuts do not intersect inside the domain. 
\end{itemize}
An example of an indecomposable, but not strongly indecomposable domain is
shown in Figure \ref{notstrong}. 

\subsection{Count of domains}
In what follows any domain is assumed to belong to
 $\mathfrak D$ and to have  Maslov index one, i.e.
 our domains have no corners with
negative multiplicities or with multiplicities bigger than $3/4$ (since
domains from $\mathfrak D$ have no corners inside).

We say that a cut touches a boundary component $A$ at an oval $B$
if either the end point of this cut belongs to $B$ or
the cut leaves $A$ along $B$.
We define the distance between two
cuts touching $A$
to be odd, if one of them touches $A$ at
a vertical oval and another one  at a  horizontal oval; 
otherwise the distance is even.

A boundary component is called {\it special} if it is an oval with
a common corner of two generators (see Figure \ref{domain} right
for an example), otherwise the component is non--special.
Let us call an inner boundary component {\it bad}, if it does not
have obtuse corners.


\begin{thm}\label{count}
Any strongly indecomposable domain without bad components
counts for the differential.
\end{thm}


The proof is given in Section 2.4 after the detailed analysis of the
 structure of
indecomposable domains.

\vspace*{.2cm}

\noindent
{\bf Remark.}
 In the long oval complex, $\partial\lo$ coincides 
with the  differential of 
link Floer homology defined by
counting of pseudo--holomorphic discs.
The differential $\partial\sh$ at least for some
order of shortening can not be realized by such count.
 This is because,
according to the Gromov compactness theorem the count of indecomposable
domains is independent of the complex structure, i.e. each of them
 either counts or not
for any complex structure.
 In our setting, 
e.g. for the domain in Figure
\ref{notstrong} one can always choose an order of shortening
in such a way that our count differs from the analytic one. Indeed,
if this domain is obtained by first shortening the dotted  part of 
the oval labeled by 1 and then the dashed part, the domain does not count.
However, if we get it by shortening the part 2, and then the dashed one, it
counts.

\subsection{Structure of domains} Here we provide some technical
results needed for the proof of Theorem \ref{count}. 
\begin{dfn}{\rm
A boundary component
 $C_1\subset \partial D_{\x,\y}$ is called 
$y$--connected  with another component
$C_2\subset\partial D_{\x,\y}$ if for any point
$y\in C_1$ and $q\in C_2$ disjoint from the corners, there exists a unique
path without self--intersections
starting at $y$ and ending at $q$, such that

$1)$ the path goes along cuts or $\partial D_{\x,\y}$, where the arcs
of horizontal and vertical ovals alternate along the path
(an arc can consist of the union of a cut and some part of $\partial D_{\x,\y}$,
as long as they are on the same oval);

$2)$ the corners of the path (i.e. intersection points of horizontal
and vertical segments)
come alternatively
from $\x$ and $\tilde \y$, 
where  $\tilde \y$ contains $\y$ and intersection points
of cuts with $\partial D_{\x,\y}$; 

$3)$ the intersection of the path with a boundary component
is neither a point nor the whole component.

$4)$ the
first corner belongs to $\x$;

}

\end{dfn}

\begin{figure}
\mbox{\epsfysize=4cm \epsffile{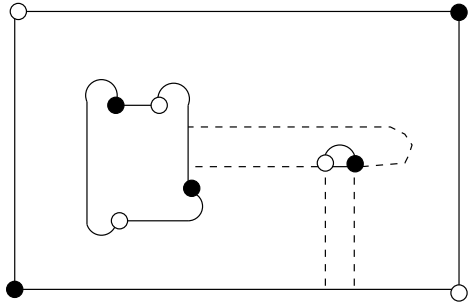}}
\v8

\caption{{\bf An indecomposable domain whose inner 
boundary components are not y--connected}}
\label{degen}
\end{figure}

As an example consider the domain shown in Figure \ref{degen}.
Let us denote the left inner boundary component by $C_1$ and the right one 
by $C_2$.
For any choice of a point $y$ on a vertical oval of  $C_1$ 
(disjoint from the corners), there is no path
$y$--connecting $C_1$ with $C_2$. On the other hand, if we choose
$y$ on the horizontal oval of $C_1$ there are two such paths. 
Hence, $C_1$ and $C_2$ are not $y$--connected.

\begin{lem}\label{new}
In  a strongly  indecomposable domain without bad components,
any two boundary components are $y$--connected.
\end{lem}

\begin{proof}
Let $c$ be the total number of boundary components in our domain.
The proof is
 by induction on  $c$. Assume first  that 
our domain has no special components. 


Further, all corners of our domain have positive multiplicities not bigger than $3/4$
and our domain has
Maslov index one. 
Since we have no bad components,  \eqref{mas-domain} implies
that every inner boundary component has exactly one obtuse corner.

Suppose  $c=2$.
If one of the cuts from the obtuse corner
connects the inner boundary component with itself, 
the domain is decomposable (see Figure \ref{02a} $(a)$).
If it is not the case, then an easy check verifies the claim (compare 
Figure \ref{02a} $(b)$).

\begin{figure}
\mbox{\epsfysize=3.5cm \epsffile{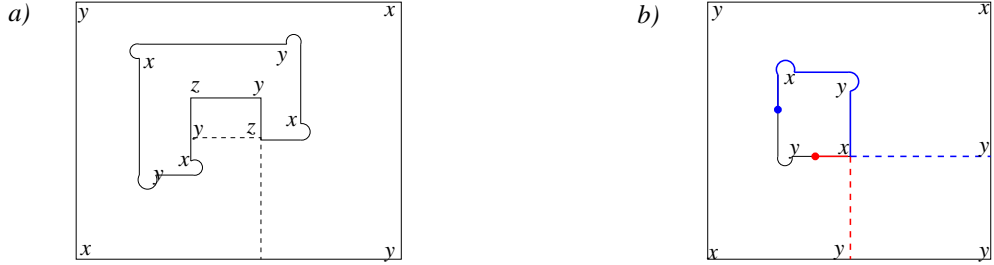}} 
\v8

\caption{{\bf Case $b=0$, $c=2$.} $a)$ The domain is the composition 
$D_{\x,\y}\circ D_{\y,\z}$, where $\z\setminus (\z\cap\y)$ are the two
 points marked by $z$.
 The cuts 
are shown
by dashed lines.  $b)$ 
Indecomposable domain. The corners 
from $\tilde\y$ are marked by $y$. For any choice of a point $y$
 on the vertical
oval of the inner component, the path $y$--connecting both components uses
a horizontal cut. Analogously, for $y$ on a horizontal oval, we have 
to use the vertical cut. Examples are shown in  blue and red, respectively. }
\label{02a}
\end{figure}

Assume the claim holds for  $c=n-1$. Suppose $c=n$,
 and our domain is
 indecomposable. 
Let us denote by $A$
the $n$--th component. Let us first assume that
  there are no cuts ending at $A$. In this case the
two cuts from the obtuse corner $y$--connect  $A$
 with some other components
which are all $y$--connected 
by  induction. 

If there is a component connected with $A$ by two cuts, then
it is $y$--connected with the outside exactly in the case when
$A$ has this property. To check this, 
 it is sufficient to find a required path  for two choices
of $y$ (before and after one corner) on this component. An example is
shown in Figure \ref{0nc}.  
Hence, when all cuts ending at $A$ come from  components 
connected  with $A$ by two cuts
(as in Figure \ref{0nc}),
 then $A$ is $y$--connected
to the outside by the previous argument.

\begin{figure}
\mbox{\epsfysize=5cm \epsffile{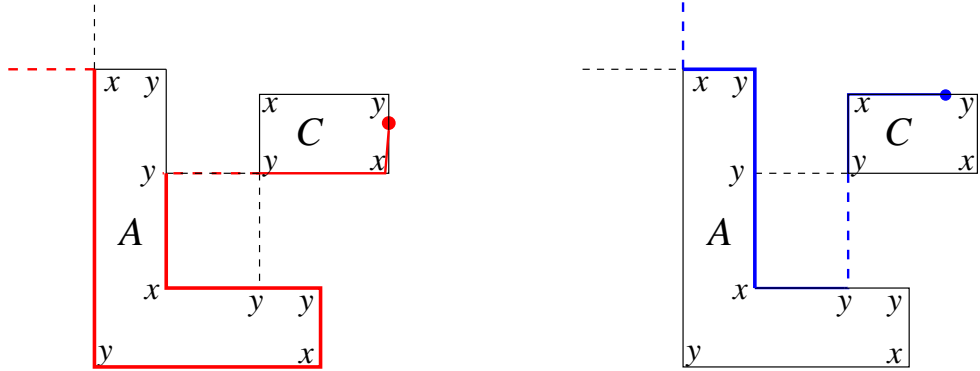}} 
\v8
\caption{{\bf $C$ is $y$--connected to $A$.} 
The two choices of $y$ are shown by red and blue dots. The connected
paths have the corresponding colors.  All corners
without cuts are assumed to be acute.}
\label{0nc}
\end{figure}


In the case, when $C\subset\partial D_{\x,\y}$ is connected with $A$ 
by just one cut or $A$ and $C$ exchange their cuts, 
 all components  are again $y$--connected 
to each other except when the following happens.
The path
  described
 in Definition 2.2
after leaving $A$ (along one of the cuts)
comes back to $A$ 
 without visiting all other components. Since this path leaves and enter
any component along cuts at odd distance (compare
Figure \ref{0na}), and has no self intersections, it  
can be used to decompose the domain, which contradicts the assumption.
\begin{figure}
\mbox{\epsfysize=3.5cm \epsffile{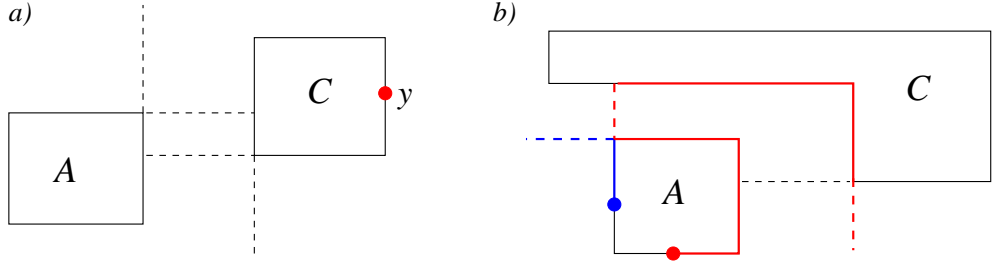}} 
\v8
\caption{{\bf  $A$ and $C$ exchange two cuts.} 
$a)$ Domain is decomposable. The point $y\in C$ is not connected with $A$
by a path described in Definition 2.2. 
$b)$
Indecomposable domain. $A$ and $C$ are  $y$--connected.
} 
\label{0na}
\end{figure}

It remains to consider the case where 
$D_{\x,\y}$ has  straight corners, or special components, i.e.
   ovals with a common corner
of two generators. 
We proceed by induction on the number of special components.
Assume that we have  one special component. 
Then the cut from the straight corner $y$--connects this
component  with any other one (otherwise not connected with  the special one) 
by the 
previous argument. 
The obvious induction  completes
the proof. 
\end{proof}

\subsection{Proof of Theorem \ref{count}}
The proof is again
 by induction on the number of boundary components $c$
in the domain. 
If  $c=1$, it is easy to see by prolonging ovals that any  
immersed polygon
counts.

\begin{figure}
\mbox{\epsfysize=5cm \epsffile{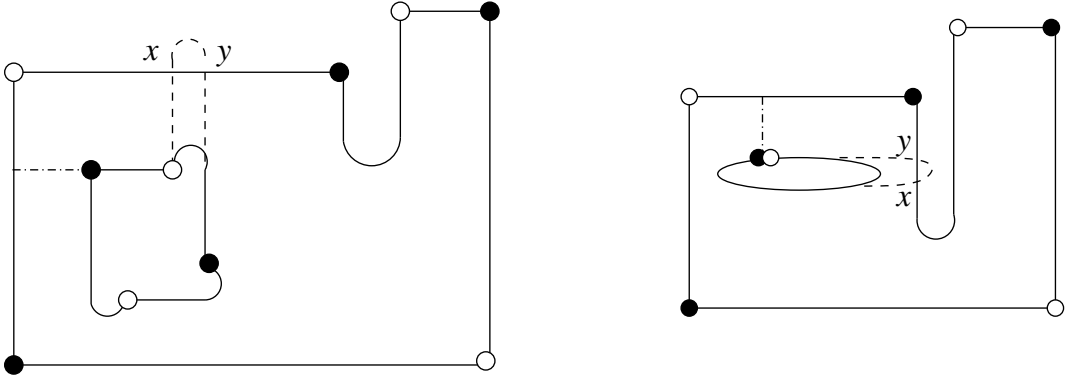}}
\v8

\caption{{\bf Strongly indecomposable domains realizing  differentials
from black to white points.} The prolongated oval is shown
by a dashed line.}
\label{domain}
\end{figure}

Assume that for $c=n-1$, the claim holds.
Suppose our complex  has a
strongly indecomposable  domain 
$D\in \frak D$  without bad components and $c=n$. 
Examples with $c=2$ are drawn in Figure \ref{domain}.

Let us stretch one oval in $D$ connecting two boundary components.
The result is 
 a domain $D'$. 
 Let  $x$ and $y$ be the corners of the bigon,
obtained after stretching.
The stretched oval connects $y$ to some  boundary
component, say $A$. By Lemma \ref{new}, 
$y$ can also be connected with $A$
by a unique path inside $D$. 
This path is not affected by the prolongation, since otherwise
 the domain would  not be  strongly
indecomposable (compare Figure \ref{notstrong}).


Hence, $D'$ can be represented
as a union of two domains
connecting some generators and having less boundary components.
The unique path connecting $y$ with $A$ leaves any 
boundary component  along a cut.
Moreover, the path has no self--intersections.
Therefore,
 $D'$ is a union of two strongly indecomposable Maslov index one domains
without bad components and with positive multiplicities at 
the corners not bigger than $3/4$. These both domains
 count for the differential 
by the induction hypothesis.
We conclude 
that the domain $D$ also counts for the differential.




\qed

\section{Computations}
In this section we 
show how $\HFK$ of small knots can be computed by hand and discuss
the computer program written by Droz.

\subsection{$5_2$ knot} 
\begin{figure}
\mbox{\epsfysize=7cm \epsffile{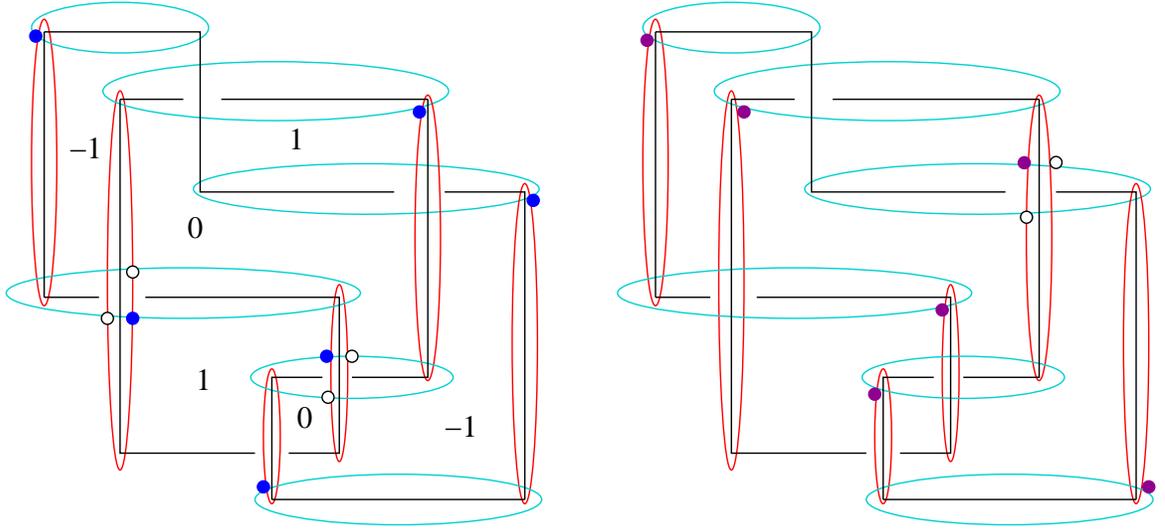}}
\v8

\caption{{\bf A complex for $5_2$ knot.} The colored dots show 
generators in the maximal Alexander grading equal to 1. A number assigned
to a region is minus the winding number for its points. }
\label{newdiag}
\end{figure}

Figure \ref{newdiag} shows
a rectangular diagram for $5_2$ knot of complexity $n=7$
obtained from the original diagram in Figure \ref{diag}
by cyclic permutations (compare \cite{Dy}).
An  advantage of this diagram is that there are no regions
counted for the differential.

The Alexander grading of a generator
is given by the formula $A(\x)=\sum_{x\in \x} a(x) -2$.
 The maximal Alexander grading is equal to one. There  are two generators
in this grading shown by colored dots in Figure \ref{newdiag}.
Both of them have Maslov grading $2$.

The homology of our complex is $\HFK (5_2)\otimes V^{6}$. Hence
in Alexander grading zero, we have 12 additional generators
coming from the multiplication with $V$.
Note that our complex  has  15 generators in Alexander grading zero.
Indeed, 12 of them can be obtained by 
moving one point of a generator in Alexander grading one
 to the other side of the oval.
In three cases, depicted by white dots there are two possibilities
to move a point. 
This gives 3 additional  generators. 
 Note that these moves
 drop Maslov index by one.
To compute $\HFK$ in the negative Alexander gradings
we  use the symmetry (\ref{sym}).

Finally, we derive that
$\HFK(5_2)$ has rank two in the  Alexander--Maslov bigrading
$(1,2)$, rank three in 
$(0,1)$,  and rank two in the bigrading
$(-1,0)$. To compare, the Alexander polynomial is 
$\Delta_{5_2}(t)=2(t+t^{-1})-3$. The knot $5_2$ is not fibered and
its Seifert genus is one.

\subsection{Droz's program} 
\begin{figure}
\mbox{\epsfysize=8cm \epsffile{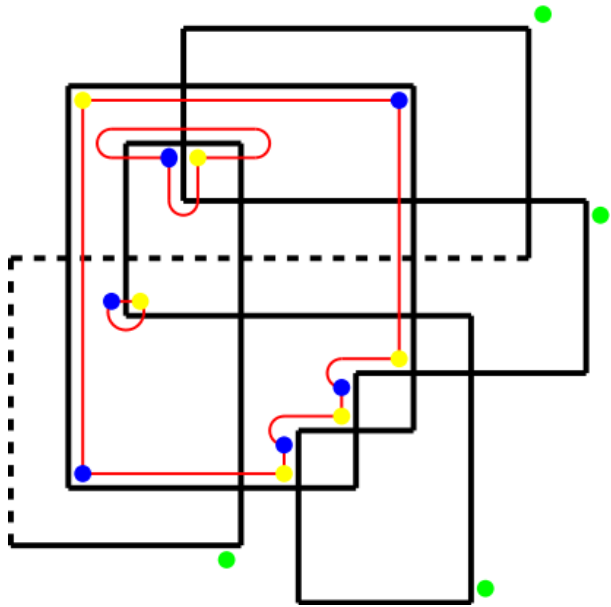}}
\v8

\caption{{\bf Domain counted twice.} The knot $11n3$ is drawn in black, 
the dashed
lines have no ovals.
The boundary of the domain from the blue to the yellow points is  red.
The common points of these two generators are green. 
If the domain is obtained
by shortening of the smallest inner component, it counts
twice over $\Z$. }
\end{figure}

In \cite{D1},
 Droz extend our construction over $\Z$ and
wrote a computer program calculating
the homology of the resulting complex \cite{Droz}. The program  is
 installed on the Bar--Natan's
Knot Atlas. As a byproduct, his program 
generates  rectangular diagrams of  knots  and links
 and  allows to change them by Cromwell--Dynnikov moves.
The program can be used to determine Seifert genus and fiberedness
of knots  until 16 crossings. 

According to Droz's computations,
  the number of generators in our complex 
is significantly smaller than that in the MOS complex.
Moreover, for small knots, almost 
all domains suitable for the differential
are embedded polygons, so they always count for the differential.
For example, for   knots admitting 
rectangular diagrams
of complexity 10,
the number of generators in the MOS complex 
 is $10!=3'628'800$. 
Our complex has on average about $50'000$ generators among them about
$1'000$ in
the  positive Alexander gradings. 
The knot $12n2000$ admits a rectangular diagram of complexity 12,
where $12!=479'001'600$.
Our complex has  $1'411'072$ generators with
 $16'065$ of them  in the positive Alexander gradings.

Furthermore, Droz's program produced 
examples of domains counted  more than once
over $\Z$. We do not know similar examples in
 the analytic setting.
In  Figure 15 one such domain is shown.
This domain has a degenerate system of cuts and its count depends
on the order of shortening of ovals. One specific order
gives multiplicity 2 for this domain.


\end{document}